\documentclass{amsart}

\setbox0=\hbox{$+$}
\newdimen\plusheight
\plusheight=\ht0
\def\+{\;\lower\plusheight\hbox{$+$}\;}

\setbox0=\hbox{$-$}
\newdimen\minusheight
\minusheight=\ht0
\def\-{\;\lower\minusheight\hbox{$-$}\;}

\setbox0=\hbox{$\cdots$}
\newdimen\cdotsheight
\cdotsheight=\plusheight%\ht0
\def\cds{\lower\cdotsheight\hbox{$\cdots$}}

\newcommand{\df}{\dfrac}

 \renewcommand{\a}{\alpha}

\renewcommand{\t}{\varphi}

\renewcommand{\(}{\left\(}
\renewcommand{\)}{\right\)}

\renewcommand{\i}{\infty}
\renewcommand{\b}{\beta}
\renewcommand{\pmod}[1]{\,(\textup{mod}\,#1)}

\newcommand\mbgr{\mbox{BG-rank}}
\newcommand{\beqs}{\begin{equation*}}
\newcommand{\eeqs}{\end{equation*}}
\numberwithin{equation}{section}
 \theoremstyle{plain}
\newtheorem{theorem}{Theorem}[section]

\newtheorem{corollary}[theorem]{Corollary}

\begin{document}

\title[\tiny{On the representations of integers by the sextenary quadratic form...}]
%$x^2+y^2+z^2+ 7s^2+7t^2+ 7u^2$ and $7$-cores}]
{On the representations of integers by the sextenary quadratic form
$x^2+y^2+z^2+ 7s^2+7t^2+ 7u^2$ and $7$-cores}
\author{Alexander Berkovich}
\author{Hamza Yesilyurt}

\address{Department of Mathematics, University of Florida, 358 Little Hall,  Gainesville, FL
  32611, USA}\email{alexb@math.ufl.edu}
\address{Department of Mathematics, Bilkent University, 06800, Bilkent/Ankara, Turkey}
\email{hamza@fen.bilkent.edu.tr}
\thanks{Research was supported in part by NSA grant H98230-07-01-0011.}

\begin{abstract}
In this paper we derive an explicit formula for the number of representations
of  an integer  by the sextenary form  $x^2+y^2+z^2+ 7s^2+7t^2+ 7u^2$.
We  establish the following intriguing  inequalities
\begin{equation*}
2b(n)\geq a_7(n) \geq b(n)\;\;\text{for}\;\;  n  \neq 0,2,6,16.
\end{equation*}
Here  $a_7(n)$  is the number of partitions of $n$ that are  7-cores and  $b(n)$ is the number of  representations of $n+2$
 by the sextenary form  $( x ^2+  y ^2+z ^2+ 7s ^2 + 7t ^2+  7u^2)/8$
with  $x,\,y,\,z,\,s,\,t$ and $u$ being odd.
\end{abstract}
\maketitle

\section{Introduction}
Recall that a partition is called  a $t$-core if it
has no rim hooks of length $t$ \cite{JK}.
Let $a_t(n)$ be the number of $t$-core partitions of $n$. It is well
known that \cite{K}, \cite{GKS}
\begin{equation}
\sum_{n \geq 0}a_t(n)q^n=\sum_{\substack{\overrightarrow{n}
 \in \mathbb{Z}^t,\;\; \overrightarrow{n}.\overrightarrow{1_t}=0}}
 q^{\tfrac{t}{2}\|\overrightarrow{n}\|^2+\overrightarrow{b_t}.\overrightarrow{n}}=\df{E^t(q^t)}{E(q)},
\end{equation}
where
\begin{equation}
\overrightarrow{b_t}:=(0,1,2,...,t-1),\quad
\overrightarrow{1_t}:=(1,1,...,1),
\end{equation}
\begin{equation*}
E(q):=\prod_{n=1}^\i(1-q^n). \notag\\
\end{equation*}

Let

\begin{equation}\label{22i}
\varphi(q) := \sum_{n=-\i}^{\i}q^{n^2},\;\;\;\psi(q)
:=\sum_{n=0}^{\i}q^{n(n+1)/2}.
\end{equation}
Throughout the paper we assume that $q$ is a complex number with $|q|<1$. For convenience, the coefficient  of $q^n$ in the expansion of $H(q)$
will be denoted as $[q^n]H(q)$. For a partition  $\pi$, BG-rank$(\pi)$ is defined as an
alternating sum of parities of parts of $\pi$ \cite{BG1}, \cite{BG2}. In
\cite{AH}, the authors found positive  $eta$-quotient
representations for the $7$-core generating functions $\sum_{n\geq 0}a_{7,j}(n)q^n$,
where $a_{7,j}(n)$ denotes the number of $7$-cores of $n$ with
$\mbgr=j$ and  established a number of
inequalities for $a_{7,j}(n)$ with $j=-1,0,1,2$ and $a_7(n)$. In this paper, we prove lower and upper bounds for $a_7(n)$, namely,

\begin{theorem}\label{TH1}
\begin{equation}\label{ubs}
[q^n]\Bigr(1+5q^6+q^{16}+2q\psi^3(q)\psi^3(q^7)\Bigl)\geq[q^n]\df{E^7(q^7)}{E(q)},
\end{equation}
where the inequality is strict if $n \neq 0,6 $ or $16$.
\end{theorem}
and
\begin{theorem}\label{TH2}
\begin{equation}\label{lbs}
[q^n]\left(\df{E^7(q^7)}{E(q)}+q^2\right)\geq[q^n]\left(q\psi^3(q)\psi^3(q^7)\right).
\end{equation}
\end{theorem}
We should remark that the inequality in \eqref{lbs} is only strict as it can be seen in the proof of Theorem \ref{TH2} when $n$ is even and $n\neq 2$.
Theorem \ref{TH1} and Theorem \ref{TH2} are proved in sections \ref{tTH1} and \ref{tTH2}.

It is well known that every integer can be written as sum of three triangular numbers, that is $[q^n]\psi^3(q)>0$  for all $n \geq 0$. This together with \eqref{lbs} implies that
\begin{equation*}
[q^n]\df{E^7(q^7)}{E(q)}>0 \;\;\text{for all}\; n \geq 3, \;\;\text{and hence for all}\;\; n \geq 0.
\end{equation*}
In fact,  Granville and Ono  showed that \cite{GO}  if  $t \geq 4$, then
\begin{equation*}
 [q^n]\df{E^t(q^t)}{E(q)}>0\;\;\text{ for all}\;\; n \geq 0.
 \end{equation*}
The lower bound given by \eqref{lbs} improves the Granville-Ono result  when $t=7$.

Essential to our proofs are the following theta function identities which we prove in section \ref{tTH3} by employing the theory of modular equations
\begin{theorem}\label{TH3}
\begin{align}
7\t^3(-q)\t^3(-q^7)
=& -49\Bigr(q^2\df{E^7(q^7)}{E(q)}+qE^3(q)E^3(q^7)\Bigl)\notag\\
&+56\Bigr(7q^4\df{E^7(q^{14})}{E(q^2)}+q^2E^3(q^2)E^3(q^{14})\Bigl)
-\df{E^7(q)}{E(q^7)}+8\df{E^7(q^2)}{E(q^{14})},\label{ntp1}
\end{align}
and
\begin{align}
56q^3\psi^3(q)\psi^3(q^7)
=& 49q^2\df{E^7(q^7)}{E(q)}+7qE^3(q)E^3(q^7)\notag\\
&-49\Bigr(q^4\df{E^7(q^{14})}{E(q^2)}+q^2E^3(q^2)E^3(q^{14})\Bigl)
+\df{E^7(q)}{E(q^7)}-\df{E^7(q^2)}{E(q^{14})}.\label{ntp2}
\end{align}
\end{theorem}
The rest of the paper is organized as follows. In the next section, we recall two Lambert series identities of Ramanujan which we extensively use in our proofs. In section \ref{med},
we give a brief introduction to modular equations. Then, we prove  Theorem \ref{TH3} and from it we derive an explicit formulas for the number of representations
of  an integer  by the sextenary forms  $x^2+y^2+z^2+ 7s^2+7t^2+ 7u^2$ and  $( x ^2+  y ^2+z ^2+ 7s ^2 + 7t ^2+  7u^2)/8$
with  $x,\,y,\,z,\,s,\,t$ and $u$ being odd for the later case. In the last two sections, Theorem \ref{TH1} and Theorem \ref{TH2} are proven.

\section{Two Lambert Series Identities of Ramanujan}
We start with two Lambert series identities of Ramanujan\cite{chan} which we
will employ in our proofs.
\begin{align}\label{Ldef}
L(q)&:=\df{8}{7}\Bigr(1-\df{E^7(q)}{E(q^7)}\Bigl)
-7qE^3(q)E^3(q^7)\notag\\
&=\sum_{n=1}^{\i}\Bigr(\df{n}{7}\Bigl)\df{n^2q^n}{1-q^n}
\end{align}
and
\begin{align}\label{Kdef}
K(q)&:=8q^2\df{E^7(q^7)}{E(q)}
+qE^3(q)E^3(q^7)\notag\\
&=\sum_{n=1}^{\i}\Bigr(\df{n}{7}\Bigl)\df{q^n(1+q^n)}{(1-q^n)^3}.
\end{align}
We should remark that \eqref{Ldef} and \eqref{Kdef}  are equivalent under the imaginary transformation \cite{chan2}.
 It is easy to see that
\begin{equation}\label{mtp}
L(q)=\sum_{n=1}^{\i}\Bigr(\sum_{d|n}d^2\Bigr(\df{d}{7}\Bigl)\Bigl)q^n
\;\;\text{and}\;\;K(q)=\sum_{n=1}^{\i}\Bigr(\sum_{d|n}d^2\Bigr(\df{n/d}{7}\Bigl)\Bigl)q^n.
\end{equation}
The coefficients of $L(q)$ and $K(q)$ are clearly multiplicative.
The reader may wish to consult \cite{apostol} for background on
multiplicative functions, convolution of multiplicative functions
and Legendre's symbol. Using multiplicity it is easy to conclude from \eqref{mtp} that

 \begin{align}
[q^n]L[q]&= (-1)^{b} \prod_{i=1}^{r} \df{1-p_i^{2v_i+2}}{ 1-p_i^2}
\prod_{j=1}^{s}\df{(-1)^{w_j}+ q_j^{2w_j+2}}{1+q_j^2},\label{form1}\\
[q^n]K[q]&= 7^{2c} \prod_{i=1}^{r} \df{1-p_i^{2v_i+2}}{ 1-p_i^2}
\prod_{j=1}^{s}\df{(-1)^{w_j}+ q_j^{2w_j+2}}{1+q_j^2},\label{form2}
\end{align}
where $n$ has the prime factorization
\begin{equation*}
n=7^c\prod_{i=1}^{r}p_i^{v_i}\prod_{j=1}^sq_j^{w_j},
\end{equation*}
with  $p_i \equiv  1,2,4\pmod {7}$,  $q_j \equiv 3,5,6 \pmod {7}$,
and $b=\sum_{j=1}^sw_j$. We note that \eqref{form2} was stated as Lemma 1 in \cite{GKS} .

Next, let
\begin{equation}\label{Mdef}
M(q):=qE^3(q)E^3(q^7).
\end{equation}
From \cite[p.~11, Lemma 2]{GKS}, we have

\begin{equation}\label{form3}
 [q^n]M(q)=\left\{ \begin{array}{ll}
         (-7)^{c}\prod_{i=1}^{r} F(p_i,v_i)
\prod_{j=1}^{s}q_i^{w_j}& \text{if each $w_j$ is even}\\
         0 &\text{otherwise,} \end{array} \right.
\end{equation}
where the prime factorization of $n$ is defined as above and
\begin{equation}\label{Fdef}
F(p,r)  := \df{\b ^{2r+2} -\bar{\b}^{2r+2} }{ \b ^2 -\bar{\b} ^2}
\end{equation}
with
\begin{equation*}
\b= x+ \sqrt{-7}y,\;\; \bar{\b}= x- \sqrt{-7}y,
\end{equation*}
where $x$ and $y$ are the positive unique integers satisfying $p=x^2+7y^2$ provided $p\equiv  1,2,4\pmod {7}$ and $p>2$. If $p=2$, then
\begin{equation*}
   \b= (1+ \sqrt{-7})/2,\;\; \bar{\b}=  (1- \sqrt{-7})/2.
\end{equation*}

Next, We give background information on modular
equations.
\section{Modular Equations}\label{med}

For $0<k<1$, the complete elliptic integral of the first
kind $K(k)$, associated with the modulus $k$, is defined by

\begin{equation*}
K(k):=\int_0^{\pi/2}\df{d\theta}{\sqrt{1-k^2\sin^2\theta}}.
\end{equation*}

The number $ k^\prime := \sqrt{1-k^2}$ is called the
\textit{complementary modulus}. Let $ K, K^\prime, L, $ and $
L^\prime $ denote complete elliptic integrals of the first kind
associated with the moduli $ k, $ $ k^\prime , $ $ \ell, $ and $
\ell^\prime , $ respectively. Suppose that
\begin{equation}\label{i2}
 n\df{K^\prime}{K} = \df{L^\prime}{L}
\end{equation}
for some positive rational integer $ n. $ A relation between $ k $
and $ \ell $ induced by \eqref{i2} is called a {\it modular equation
of degree $ n.$ } There are several definitions of a modular
equation in the literature. For example, see the books by
R.~A.~Rankin \cite[p.~76]{rankin} and B.~Schoeneberg
\cite[pp.~141--142]{sc}. Following Ramanujan, set
$$ \a = k^2 \qquad \text{and} \qquad  \b = \ell^2. $$
We often say that $ \b $ has degree $ n $ over $ \a. $ If
\begin{equation}\label{i3}
q = \exp(-\pi K^{\prime}/K),
\end{equation}
two of the most fundamental relations in the theory of elliptic
functions are given by the formulas \cite[pp.~101--102]{III},
\begin{equation}\label{i4}
 \t^2(q) =\df{2}{\pi}K(k)
 \;\;\text{and}\;\;\a=k^2=1-\df{\t^4(-q)}{\t^4(q)}.
\end{equation}
The equation \eqref{i4} and elementary theta function identities
 make it possible
to write each modular equation as a theta function identity.
Ramanujan derived an extensive ``catalogue'' of formulas
\cite[pp.~122--124]{III} giving  the ``evaluations'' of $ E(q)$,
$\t(q)$, $\psi(q)$, and $ \chi(q)$ at various powers of the
arguments in terms of
$$ z := z_1 :=\df{2}{\pi}K(k) , \quad \a, \quad
\text{and}\quad q.$$ The evaluations that will be needed in this
paper are as follows
\begin{align}
\phi(-q)&=\sqrt{z}\{(1-\a)\}^{1/4},\label{list444}\\
\psi(-q)&=q^{-1/8}\sqrt{\tfrac{1}{2}z}\{\a(1-\a)\}^{1/8},\label{list5}\\
E(-q)&=2^{-1/6}q^{-1/24}\sqrt{z}\{\a(1-\a)\}^{1/24},\label{list65}\\
E(q^2)&=2^{-1/3}q^{-1/12}\sqrt{z}\{\a(1-\a)\}^{1/12},\label{list7}\\
E(q^4)&=4^{-1/3}q^{-1/6}\sqrt{z}{\a}^{1/6}{(1-\a)}^{1/24}.\label{list777}
\end{align}
We should remark that in the notation of \cite{III}, $E(q)=f(-q)$.
If $ q $ is replaced by $ q^n$, then the evaluations are given in
terms of
 $$  z_n := \df{2}{\pi}K(l), \quad \b, \quad
\text{and} \quad q^n,$$ where $ \b $ has degree $ n $ over $ \a$.

Lastly, the multiplier $ m $ of degree $ n $ is defined by
\begin{equation}\label{i5}
 m =\df{\t^2(q)}{\t^2(q^n)}=\df{z}{z_n}.
\end{equation}
The proofs of the following  modular equations of degree 7 can be
found in \cite[p.~314, Entry 19(i),(iii), (viii)]{III}

\begin{align}
(\a \b)^{1/8}+\bigr\{(1-\a)(1-\b)\bigl\}^{1/8}&=1,\label{71}\\
m=\df{1-4\Bigr(\df{\b^7(1-\b)^7}{\a(1-\a)}\Bigl)^{1/24}}{\bigr\{(1-\a)(1-\b)\bigl\}^{\tfrac{1}{8}}-(\a
\b)^{\tfrac{1}{8}}},\;\;\;&\df{7}{m}=-\df{1-4\Bigr(\df{\a^7(1-\a)^7}{\b(1-\b)}\Bigl)^{1/24}}{\bigr\{(1-\a)(1-\b)\bigl\}^{\tfrac{1}{8}}-(\a
\b)^{\tfrac{1}{8}}},\label{74}\\
m-7/m=2\bigr((\a
\b)^{1/8}-\bigr\{(1-\a)(1-\b)\bigl\}^{1/8}\bigl)&\bigr(2+(\a
\b)^{1/4}+\bigr\{(1-\a)(1-\b)\bigl\}^{1/4}\bigl).\label{775}
\end{align}

\section{Proof of Theorem \ref{TH3}}\label{tTH3}
In the language of modular equations the identities  \eqref{ntp1} and   \eqref{ntp2} are reciprocals of each other \cite[p.~216, Entry 24(v)]{III} and so we only prove \eqref{ntp2}. In
\eqref{ntp2}, we replace $q$ by $-q$ and use the evaluations given
in \eqref{list5}--\eqref{list7}, we find that

\begin{align}
-7&\sqrt{z^3}\sqrt{z_7^3}\bigr\{\a\b(1-\a)(1-\b)\bigl\}^{3/8}\notag\\
=&\df{49}{2}\df{\sqrt{z_7^7}}{\sqrt{z}}\Bigr(\df{\b^7(1-\b)^7}{\a(1-\a)}\Bigl)^{1/24}
-\df{7}{2}\sqrt{z^3}\sqrt{z_7^3}\bigr\{\a\b(1-\a)(1-\b)\bigl\}^{1/8}\notag\\
&-\df{49}{4}\df{\sqrt{z_7^7}}{\sqrt{z}}\Bigr(\df{\b^7(1-\b)^7}{\a(1-\a)}\Bigl)^{1/12}
-\df{49}{4}\sqrt{z^3}\sqrt{z_7^3}\bigr\{\a\b(1-\a)(1-\b)\bigl\}^{1/4}\notag\\
&+\df{1}{2}\df{\sqrt{z^7}}{\sqrt{z_7}}\Bigr(\df{\a^7(1-\a)^7}{\b(1-\b)}\Bigl)^{1/24}-
\df{1}{4}\df{\sqrt{z^7}}{\sqrt{z_7}}\Bigr(\df{\a^7(1-\a)^7}{\b(1-\b)}\Bigl)^{1/12}.\label{zorla}
\end{align}
We divide both sides of \eqref{zorla} by $\sqrt{z^3}\sqrt{z_7^3}$
and use \eqref{i5} and conclude that \eqref{zorla} is equivalent to
\begin{align}
&\df{49}{m^2}\Bigr\{1-\Bigr(1-\Bigr(\df{\b^7(1-\b)^7}{\a(1-\a)}\Bigl)^{1/24}\Bigl)^2\Bigl\}+
m^2\Bigr\{1-\Bigr(1-\Bigr(\df{\a^7(1-\a)^7}{\b(1-\b)}\Bigl)^{1/24}\Bigl)^2\Bigl\}\notag\\
&+7\Bigr\{4\bigr\{\a\b(1-\a)(1-\b)\bigl\}^{3/8}-2\bigr\{\a\b(1-\a)(1-\b)\bigl\}^{1/8}-7\bigr\{\a\b(1-\a)(1-\b)\bigl\}^{1/4}\Bigl\}=0.\label{zorla2}
\end{align}
We prove \eqref{zorla2}.

Set $t:=(\a \b)^{1/8}$. Then, by \eqref{71}, we have
\begin{equation}
\bigr\{(1-\a)(1-\b)\bigl\}^{1/8}=1-t.
\end{equation}
Let $x:=1-2t$, from \eqref{74}, we have
\begin{equation}\label{776}
\Bigr(\df{\b^7(1-\b)^7}{\a(1-\a)}\Bigl)^{1/24}=\df{1-xm}{4}\;\;
\text{and}\;\;\Bigr(\df{\a^7(1-\a)^7}{\b(1-\b)}\Bigl)^{1/24}=\df{1+7x/m}{4}.
\end{equation}
Similarly, \eqref{775} is equivalent to
\begin{equation}\label{777}
m-7/m=2(2t-1)(2t^2-2t+3).
\end{equation}

Now using \eqref{776} and \eqref{777},  we find after some algebra
that
\begin{align}
&\df{49}{m^2}\Bigr\{1-\Bigr(1-\Bigr(\df{\b^7(1-\b)^7}{\a(1-\a)}\Bigl)^{1/24}\Bigl)^2\Bigl\}+
m^2\Bigr\{1-\Bigr(1-\Bigr(\df{\a^7(1-\a)^7}{\b(1-\b)}\Bigl)^{1/24}\Bigl)^2\Bigl\}\notag\\&=\df{7}{16}\Bigr((m-7/m)^2+6x(m-7/m)-14x^2+14\Bigl)\notag\\
&=7\bigr(4\,{t}^{6}-12\,{t}^{5}+19\,{t}^{4}-18\,{t}^{3}+5\,{t}^{2}+2\,t\bigl).
\end{align}
Moreover,
\begin{align}
&7\Bigr\{4\bigr\{\a\b(1-\a)(1-\b)\bigl\}^{3/8}-2\bigr\{\a\b(1-\a)(1-\b)\bigl\}^{1/8}-7\bigr\{\a\b(1-\a)(1-\b)\bigl\}^{1/4}\Bigl\}\notag\\
&=7\bigr(4t^3(1-t)^3-2t(1-t)-7t^2(1-t)^2\bigl)\notag\\&=-7\bigr(4\,{t}^{6}-12\,{t}^{5}+19\,{t}^{4}-18\,{t}^{3}+5\,{t}^{2}+2\,t\bigl).
\end{align}
This completes the proof of \eqref{zorla2}. Hence the proof of
Theorem \ref{TH3} is complete. Next, we determine $[q^n]\bigr(\t^3(q)\t^3(q^7)\bigl)$ and $[q^n]\bigr(q^3\psi^3(q)\psi^3(q^7)\bigl)$  which are  the number of representations of  an integer  by the sextenary forms   $x^2+y^2+z^2+ 7s^2+7t^2+ 7u^2$ and  $( x ^2+  y ^2+z ^2+ 7s ^2 + 7t ^2+  7u^2)/8$
with  $x,\,y,\,z,\,s,\,t$ and $u$ being odd for the later case.

\begin{corollary}\label{repphi}
Suppose $n$ has the prime factorization
\begin{equation*}
n=7^c2^d\prod_{i=1}^{r}p_i^{v_i}\prod_{j=1}^sq_j^{w_j},
\end{equation*}
with  $p_i$ odd  $p_i \equiv  1,2,4\pmod {7}$,  $q_j \equiv 3,5,6 \pmod {7}$,
and $b=\sum_{j=1}^sw_j$.\\
\smallskip

\noindent if $n$ is odd, then
\begin{align}
[q^n]\left(\t^3(q)\t^3(q^7)\right)&=\df{1}{8}\bigr(7^{2c+1}- (-1)^{b}\bigl) \prod_{i=1}^{r} \df{1-p_i^{2v_i+2}}{ 1-p_i^2}
\prod_{j=1}^{s}\df{(-1)^{w_j}+ q_j^{2w_j+2}}{1+q_j^2}\notag\\&\;\;+\df{21}{4}(-7)^c\prod_{i=1}^{r} F(p_i,v_i)
\prod_{j=1}^{s}\df{q_i^{w_j}(1+(-1)^{w_j})}{2},\label{crd1}
\end{align}
if $n$ is even, then
\begin{align}
[q^n]\left(\t^3(q)\t^3(q^7)\right)&=\df{1}{24}\bigr(7^{2c+1}- (-1)^{b}\bigl)\left(4^{d+1}-7\right) \prod_{i=1}^{r} \df{1-p_i^{2v_i+2}}{ 1-p_i^2}
\prod_{j=1}^{s}\df{(-1)^{w_j}+ q_j^{2w_j+2}}{1+q_j^2}\notag\\&\;\;-\df{3}{4}(-7)^c\left(7F(2,d)+8F(2,d-1)\right)\prod_{i=1}^{r} F(p_i,v_i)
\prod_{j=1}^{s}\df{q_i^{w_j}(1+(-1)^{w_j})}{2}.\label{crd2}
\end{align}
\end{corollary}
\begin{proof}
From \eqref{ntp1} with $q$ replaced by $-q$, and the definitions \eqref{Ldef}, \eqref{Kdef}, and \eqref{Mdef}, we have that
\begin{equation}
8\t^3(q)\t^3(q^7)=8+L(-q)-7K(-q) -8L(q^2)+56K(q^2)-42M(-q)
 -48M(q^2).\label{btyi}
\end{equation}
Therefore,
\begin{equation}\label{SdS1}
8[q^{2n+1}]\bigr(\t^3(q)\t^3(q^7)\bigl)=[q^{2n+1}]\bigr(7K(q)-L(q)+42M(q)\bigl)
\end{equation}
and
\begin{equation}\label{SdS2}
8[q^{2n}]\bigr(\t^3(q)\t^3(q^7)\bigl)=[q^{2n}]\bigr(8+L(q)-8L(q^2)-7K(q)+56K(q^2)-42M(q)-48M(q^2)\bigl).
\end{equation}
These two equations together with \eqref{form1}, \eqref{form2} and \eqref{form3} imply \eqref{crd1} and \eqref{crd2}.

\end{proof}
\begin{corollary}\label{reppsi}
Suppose $n$ has the prime factorization
\begin{equation*}
n=7^c2^d\prod_{i=1}^{r}p_i^{v_i}\prod_{j=1}^sq_j^{w_j},
\end{equation*}
with  $p_i$ odd  $p_i \equiv  1,2,4\pmod {7}$,  $q_j \equiv 3,5,6 \pmod {7}$,
and $b=\sum_{j=1}^sw_j$.\\
\smallskip

\noindent if $n$ is odd, then
\begin{align}
[q^n]\left(q^3\psi^3(q)\psi^3(q^7)\right)&=\df{1}{64}\bigr(7^{2c+1}- (-1)^{b}\bigl) \prod_{i=1}^{r} \df{1-p_i^{2v_i+2}}{ 1-p_i^2}
\prod_{j=1}^{s}\df{(-1)^{w_j}+ q_j^{2w_j+2}}{1+q_j^2}\notag\\&\;\;-\df{3}{32}(-7)^c\prod_{i=1}^{r} F(p_i,v_i)
\prod_{j=1}^{s}\df{q_i^{w_j}(1+(-1)^{w_j})}{2},\label{crd1m}
\end{align}
if $n$ is even, then
\begin{align}
[q^n]\left(q^3\psi^3(q)\psi^3(q^7)\right)&=\df{1}{64}4^d\bigr(7^{2c+1}- (-1)^{b}\bigl) \prod_{i=1}^{r} \df{1-p_i^{2v_i+2}}{ 1-p_i^2}
\prod_{j=1}^{s}\df{(-1)^{w_j}+ q_j^{2w_j+2}}{1+q_j^2}\notag\\&\;\;-\df{3}{32}(-7)^c\left(F(2,d)+7F(2,d-1)\right)\prod_{i=1}^{r} F(p_i,v_i)
\prod_{j=1}^{s}\df{q_i^{w_j}(1+(-1)^{w_j})}{2}.\label{crd2m}
\end{align}
\end{corollary}
The proof of Corollary \ref{reppsi} is very similar to that of Corollary \ref{repphi} and we forgo its proof.
\section{Proof of Theorem \ref{TH1}}\label{tTH1}
From \eqref{ntp2}  and the definitions \eqref{Ldef}, \eqref{Kdef}, and \eqref{Mdef}, we have that
\begin{align}
&32q^2\Bigr(2q\psi^3(q)\psi^3(q^7)-\df{E^7(q^7)}{E(q)}\Bigl)\notag\\&=3K(q)-7K(q^2)-L(q)-2M(q)+L(q^2)-42M(q^2).\label{epl}
\end{align}
Explicit check shows that \eqref{ubs} is valid for $n=0,\,6$ or $n=16$ and so we assume in this section that $n\neq 0,\,6$ or 16. From \eqref{epl}, we see that
\begin{equation}
32[q^{2n-1}]\Bigr(2q\psi^3(q)\psi^3(q^7)-\df{E^7(q^7)}{E(q)}\Bigl)=[q^{2n+1}]\left(3K(q)-L(q)-2M(q)\right).\label{nib}
\end{equation}
Let $r(n):=[q^n]\left(3K(q)-L(q) -2M(q)\right)$.
Instead of proving that \eqref{nib} is nonnegative, we will prove the stronger statement that if $n>1$ , then
\begin{equation}
r(n)> 0.
\end{equation}

If $[q^n]M(q)=0$, then by \eqref{form1} and \eqref{form2}, we have that
 \begin{align}
r(n)&=[q^n]\Bigr(3K(q)-L(q) \Bigl)\notag\\&=\bigr(3.7^{2c}- (-1)^{b}\bigl) \prod_{i=1}^{r} \df{1-p_i^{2v_i+2}}{ 1-p_i^2}
\prod_{j=1}^{s}\df{(-1)^{w_j}+ q_j^{2w_j+2}}{1+q_j^2}>0,\label{ygh}
\end{align}
where $n$ has the prime factorization
\begin{equation*}
n=7^c\prod_{i=1}^{r}p_i^{v_i}\prod_{j=1}^sq_j^{w_j},
\end{equation*}
with  $p_i \equiv  1,2,4\pmod {7}$,  $q_j \equiv 3,5,6 \pmod {7}$,
and $b=\sum_{j=1}^sw_j$.

Let $s(n):=[q^n]M(q)$, assuming now that $s(n) \neq 0$, we have by \eqref{form1}, \eqref{form2}, \eqref{form3},

\begin{align}\label{inthis}
&r(n)=|s(n)|7^{-c}\Bigr\{\bigr(3.7^{2c}- 1\bigl) \prod_{i=1}^{r} \df{1-p_i^{2v_i+2}}{ (1-p_i^2)|F(p_i,v_i)|}
\prod_{j=1}^{s}\df{1+ q_j^{2w_j+2}}{q_j^{w_j}(1+q_j^2)}-2.7^c.\df{s(n)}{|s(n)|}\Bigl\}.
\end{align}
From \eqref{Fdef}, we observe that
\begin{equation*}
F(p,r)  = \df{\b ^{2r+2} -\bar{\b}^{2r+2} }{ \b ^2 -\bar{\b} ^2}=\b^{2r}+\b^{2r-2}\bar{\b}^2+...+\bar{\b}^{2r},
\end{equation*}
where $p=\b\bar{\b}$. Therefore,
\begin{equation}\label{fgn}
|F(p,r)|\leq (r+1)p^r.
\end{equation}
It is easy  to show that if $p$ and $q$ as above and $w$ is even, then
%\begin{align}
%&\text{if}\; p>2,\; \text{then} \df{1-p^{2r+2}}{ |F(p,r)|(1-p^2) } > 9,\label{LIST1}\\
%&\text{if}\; r>2,\; \text{then} \df{1-2^{2r+2}}{|F(2,r)|(1-2^2) } > 28,\label{LIST2}\\
%&\df{1-2^{2+2}}{|F(2,1)|(1-2^2) } =5/3,\label{LIST3}\\
%&\text{if}\;  w\; \text{is even and positive}, \;\text{then}\;\;\df{q^{2w+2}+ 1}{   q^w(1+q^2)  } > 8.\label{LIST4}
%\end{align}

\begin{equation}
\df{1-p^{2v+2}}{ |F(p,v)|(1-p^2) }\geq\left\{ \begin{array}{ll}
        11  &  \text{if}\; p>2 \;\text{or}\; v>2\\
         5/3 &\text{if}\; p=2,v=1\\
         21/5 & \text{if}\; p=2, v=2
\end{array} \right.\label{LIST1}
\end{equation}
and
\begin{equation}
\df{q^{2w+2}+ 1}{   q^w(1+q^2)  }\geq\left\{ \begin{array}{ll}
        11 & \text{if}\;q>3 \;\text{or}\; w > 3\\
         8 &\text{if}\; q=3,w=2
         .\end{array} \right.\label{LIST4}
\end{equation}

Using \eqref{LIST1} and \eqref{LIST4} in \eqref{inthis}, we conclude that
\begin{equation}
r(n)\geq 7^{-c}\left(\bigr(3.7^{2c}- 1\bigl).\df{5}{3} -2.7^c\right)>7^{-c}\left(\bigr(3.7^{2c}- 1\bigl) -2.7^c\right)\geq 0.
\end{equation}

Next, we look at even-indexed coefficients. From \eqref{epl}, we find that
\begin{align}
&32[q^{2n}]\Bigr(2q\psi^3(q)\psi^3(q^7)-\df{E^7(q^7)}{E(q)}\Bigl)\notag\\
&=[q^{2n+2}]\left(3K(q)-7K(q^2)-L(q)-2M(q)+L(q^2)-42M(q^2)\right).\label{epll}
\end{align}
Therefore, it remains to prove
\begin{equation}\label{ubs22}
[q^n]\left(3K(q)-7K(q^2)-L(q)-2M(q)+L(q^2)-42M(q^2)\right)>0,
\end{equation}
where $n$ is an even integer, $n\neq 0+2=2,\; 6+2=8\; \text{or}\; 16+2=18$. Suppose n has the prime factorization
\begin{equation*}
n=7^c2^d\prod_{i=1}^{r}p_i^{v_i}\prod_{j=1}^sq_j^{w_j},
\end{equation*}
where  $p_i$ odd $p_i \equiv  1,2,4\pmod {7}$,  $q_j \equiv 3,5,6 \pmod {7}$,
$b=\sum_{j=1}^sw_j$ and $d>0$.
Employing \eqref{form1} and \eqref{form2}, we find that
\begin{align*}
&[q^{n}]\left(3K(q)-7K(q^2)-L(q)+L(q^2)\right)\\&=\df{1}{3}\left(7^{2c}(5.4^d+4)-(-1)^b.3.4^{d}\right) \prod_{i=1}^{r} \df{1-p_i^{2v_i+2}}{ 1-p_i^2}
\prod_{j=1}^{s}\df{(-1)^{w_j}+ q_j^{2w_j+2}}{1+q_j^2}>0,
\end{align*}
which  proves \eqref{ubs22}  if $[q^n]M(q)=0$. Thus, we assume now that $[q^n]M(q) \not =0$ that is $w_i$ and hence $b$ are all even, by \eqref{form3}, we find that
\begin{align}
&[q^{n}]\left(3K(q)-7K(q^2)-L(q)-2M(q)+L(q^2)-42M(q^2)\right)\notag\\
&=\prod_{j=1}^{s}q_j^{w_j}\prod_{i=1}^{r} |F(p_i,v_i)|\Bigl\{ \df{1}{3}\left(7^{2c}(5.4^d+4)-3.4^{d}\right)\prod_{i=1}^{r} \df{1-p_i^{2v_i+2}}{ (1-p_i^2)|F(p_i,v_i)|}
\prod_{j=1}^{s}\df{1+ q_j^{2w_j+2}}{(1+q_j^2)q_j^{w_j}}\notag\\&\;\;-2.(-7)^{c}\bigl(F(2,d)+21F(2,d-1)\bigr)\prod_{i=1}^{r}\df{ F(p_i,v_i)} {|F(p_i,v_i)|} \Bigr\}\notag\\&\geq \df{1}{3}\left(7^{2c}(5.4^d+4)-3.4^{d}\right)\prod_{i=1}^{r} \df{1-p_i^{2v_i+2}}{ (1-p_i^2)|F(p_i,v_i)|}
\prod_{j=1}^{s}\df{1+ q_j^{2w_j+2}}{(1+q_j^2)q_j^{w_j}}\notag\\&\;\;-2.7^{c}\bigl|F(2,d)+21F(2,d-1)\bigr|.\label{bis}
\end{align}
Let
\begin{equation}
S_1:=\df{1}{3}\left(7^{2c}(5.4^d+4)-3.4^{d}\right)A(n)-2.7^c\bigl|F(2,d)+21F(2,d-1)\bigr|,
\end{equation}
where
\begin{equation}
A(n):=\prod_{i=1}^{r} \df{1-p_i^{2v_i+2}}{ (1-p_i^2)|F(p_i,v_i)|}
\prod_{j=1}^{s}\df{1+ q_j^{2w_j+2}}{(1+q_j^2)q_j^{w_j}}.
\end{equation}

From \eqref{LIST1}, \eqref{LIST4}, and \eqref{fgn}, we find that

\begin{align}
&S_1\geq S_2:=\df{1}{3}\left(7^{2c}(5.4^d+4)-3.4^{d}\right)-2.7^{c}\bigl|F(2,d)+21F(2,d-1)\bigr|\\
&>S_3:=\df{1}{3}\left(7^{2c}(5.4^d+4)-3.4^{d}\right)-2.7^{c}((d+1)2^d+21d2^{d-1})\\
&=\df{7^c}{3}\left(7^{c}(5.4^d+4)-3.4^{d}-3(23d+2)2^{d}\right).
\end{align}
It is easy to show that $S_3>0$ if $c \geq 1$ except for $c=1$ and $d=2$ but $S_2>0$ for $c=1$ and $d=2$. Observe that $S_3>0$ if $d>8$ and $c=0$. Direct evaluation shows that $S_2>0$ if $c=0$, $d=4,5,6,7$ or $8$. For the remaining cases, $c=0$, $d=1$,2, or 3, by \eqref{LIST1} and \eqref{LIST4}, we have that $A(n)\geq 11$ unless $n=2^d3^w$ with $0 \leq w \leq 2$. The validity of \eqref{ubs22}  can easily be checked for $n=  4,\,12,\,36$ ($n=2$,\,8 or 18 and $w=1$ are already excluded). Assuming $A(n)\geq 11$, direct computation show that $S_1>0$ if  $c=0$, $d=$1,2, or 3.  Hence, the proof of Theorem \ref{TH1} is complete.

\section{Proof of Theorem \ref{TH2}}\label{tTH2}

From \eqref{ntp2}, and the definitions \eqref{Ldef}, \eqref{Kdef}, and \eqref{Mdef}, we find that
\begin{align}
&64q^2\left(\df{E^7(q^7)}{E(q)}-q\psi^3(q)\psi^3(q^7)\right)\notag\\
&=L(q)-L(q^2)+K(q)-K(q^2)-2M(q)+6\left(K(q^2)+7M(q^2)\right)+2K(q^2).
\end{align}
It is clear from \eqref{form2} that $[q^n]K(q^2)>0$. Below we assume that $n \not= 2$.  Validity of  \eqref{lbs} for the corresponding value of $n=2$ can easily be checked. Therefore, it suffices to prove that if $n \not =2$, then
\begin{equation}\label{xc1}
[q^n]\left(K(q)+7M(q)\right) > 0
\end{equation}
and
\begin{equation}\label{xc2}
[q^n]\left(L(q)-L(q^2)+K(q)-K(q^2)-2M(q)\right) \geq 0.
\end{equation}
We start with \eqref{xc1}. From \eqref{form2} ,
\begin{equation}
[q^n](K(q)+ 7M(q)) \geq 0 \;\; \text{if}\;\; s(n):=[q^n]M(q)=0.
\end{equation}
Assuming that $s(n) \not =0$, by \eqref{form2}, \eqref{form3}, \eqref{LIST1},  and by \eqref{LIST4}, we find that
\begin{equation}
\df{[q^n] K(q)}{|s(n) |}=7^c\prod_{i=1}^{r} \df{1-p_i^{2v_i+2}}{ (1-p_i^2)|F(p_i,v_i)|}
\prod_{j=1}^{s}\df{1+ q_j^{2w_j+2}}{q_j^{w_j}(1+q_j^2)}>7,
\end{equation}
provided $n \not = 2$ or 4 . However, $s(4)>0$ and so  we conclude that
 \begin{equation}
[q^n]\left(K(q)+7M(q)\right)=|s(n) |\left (   \df{[q^n] K(q)}{|s(n) |}   +  7\df{s(n)}{|s(n)|}\right)>0.
\end{equation}
Next, we prove \eqref{xc2}. Assume as before that $n$ has the prime factorization,

\begin{equation*}
n=7^c2^d\prod_{i=1}^{r}p_i^{v_i}\prod_{j=1}^sq_j^{w_j},
\end{equation*}
with $p_i$ odd $p_i \equiv  1,2,4\pmod {7}$,  $q_j \equiv 3,5,6 \pmod {7}$,
$b=\sum_{j=1}^sw_j$.\\
From \eqref{form1} and \eqref{form2}, we find that
\begin{equation*}
[q^n]\left(L(q)-L(q^2)+K(q)-K(q^2)\right)=(7^{2c}+(-1)^b)2^{2d}\prod_{i=1}^{r} \df{1-p_i^{2v_i+2}}{ (1-p_i^2)}
\prod_{j=1}^{s}\df{(-1)^{w_j}+ q_j^{2w_j+2}}{(1+q_j^2)} \geq 0,
\end{equation*}
which proves \eqref{xc2} if $s(n)=0$. Next assume that $s(n) \neq 0$. Then, $w_j$ and $b$ are all even and by employing \eqref{form3}, \eqref{LIST1},  and \eqref{LIST4}, we find that
\begin{align}
&[q^n]\left(L(q)-L(q^2)+K(q)-K(q^2)-2M(q)\right)\\
=&|s(n)|\Bigl\{\df{(7^c+7^{-c})2^{2d}}{|F(2,d)|}\prod_{i=1}^{r} \df{1-p_i^{2v_i+2}}{ (1-p_i^2)|F(p_i,v_i)|}
\prod_{j=1}^{s}\df{1+ q_j^{2w_j+2}}{q_j^{w_j}(1+q_j^2)}-2\df{s(n)}{|s(n)|}\Bigr\} \geq 0,
\end{align}
since  $|F(2,d)|\leq(d+1)2^d\leq 2^{2d}$ by \eqref{fgn}.

\section{Concluding Remarks}
We would like to point out that another upper bound for the coefficients of $7$-cores is given by the inequality
\begin{equation}
[q^n]\Bigr(\t^3(q)\t^3(q^7)\Bigl)>5[q^n]\left(q^3\psi^3(q)\psi^3(q^7)\right)>[q^n]\Bigr(q^2\df{E^7(q^7)}{E(q)}\Bigl)\; \text{for}\; n \neq 2,4,7,14,22,29,58.
\end{equation}
 The proof of the first part of this inequality is similar to that of Theorem \ref{TH1} and is omitted for space considerations. The second part of this inequality follows from Theorem \ref{TH1}. It would be interesting to prove all these inequalities for $7$-cores in a completely elementary manner. It is natural to ask if  our  inequalities   extend to general $t$-cores. We offer the following inequality as a conjecture:

\begin{equation*}
[q^n] \Bigr( \psi(q)\psi(q^t)\Bigl)^{(t-1)/2}\geq [q^n]\Bigr(  \df{E^t(q^t)}{E(q)}\Bigl),
\end{equation*}
valid for all $n$, provided that $t$ is an odd integer greater or equal to 11.

\section{Acknowledgment}
We would like to thank George Andrews, Frank Garvan, Michael Somos and James Sellers for their interest and helpful comments.


\begin{thebibliography}{99}

\bibitem{apostol}
T.~M.~Apostol, \emph{Introduction to Analytic Number Theory},
Springer--Verlag, New York, 1976 .

\bibitem{BG1}
A. Berkovich, F. G. Garvan, \emph{On the Andrews-Stanley refinement
of Ramanujan's congruence modulo $5$ and generalization}, Trans.
Amer. Math. Soc. \textbf{358} (2006), 703--726.

\bibitem{BG2}
A. Berkovich, F. G. Garvan, \emph{The BG-rank of a partition and its applications},
Adv. in Appl. Math. 40 (2008), no. 3, 377--400.


\bibitem{AH}
A.~Berkovich, H.~Yesilyurt, \emph{New Identities  for 7-cores with prescribed BG-rank}, to appear in Discrete
Math.

\bibitem{III}
B.~C.~Berndt, \emph{Ramanujan's Notebooks, Part III},
Springer--Verlag, New York, 1991.

\bibitem{chan}
H.~H.~Chan, \emph{New proofs of Ramanujan's partition identities for moduli 5 and 7}, J.~ Number Theory  53 (1995) ,144-159
\bibitem{chan2}
H.~H.~Chan, \emph{On the equivalence of Ramanujan's partition identities and a connection with the Rogers-Ramanujan continued fraction}, J. Math. Anal. Appl. 198 (1996), no. 1, 111--120

\bibitem{GKS}
F.~Garvan, D.~Kim and D.~Stanton, \emph{Cranks and t-cores},
Invent.~Math.~\textbf{101} (1990), 1--17.

\bibitem{GO}
A.~Granville and K.~Ono, \emph{Defect zero $p$-blocks for finite
simple groups}, Trans.~Amer.~Math.~Soc.~\textbf{348}(1) (1996),
331-347.


\bibitem{JK}
G.~James and A.~Kerber, \emph{The Representation Theory of the
Symmetric Group}, Addison-Wesley, Reading, MA, 1981


\bibitem{K}
A.~A.~Klyachko, \emph{Modular forms and representations of symmetric
groups}, Jour.~Soviet.~Math.~\textbf{26} (1984), 1879-1887.

\bibitem{rankin}
R.~A.~Rankin, \emph{Modular Forms and Functions}, Cambridge
University Press, Cambridge, 1977.


\bibitem{sc}
B.~Schoeneberg, \emph{Elliptic Modular Functions}, Springer-Verlag,
New York, 1974.

\end{thebibliography}
\end{document}